\documentclass[preprint,12pt]{elsarticle}
\usepackage{graphicx}
\usepackage{amssymb}
\usepackage{graphicx,setspace}
\usepackage{psfrag}
\usepackage{amsfonts}
\expandafter\let\csname equation*\endcsname\relax
\expandafter\let\csname endequation*\endcsname\relax
\usepackage{amsmath}
\usepackage{amssymb, amsthm}
\usepackage{mathrsfs}
\usepackage{epstopdf}
\usepackage{float}
\usepackage{lineno,hyperref}
\usepackage{color}
\usepackage{subfigure}
\usepackage{dsfont} 
\usepackage{amssymb} 
\usepackage{graphicx,color}
\usepackage{extarrows}
\usepackage{epstopdf}

\journal{Applied Mathematics and Computation}

\begin{document}

\begin{frontmatter}



\title{Rare events in a stochastic vegetation-water dynamical system based on machine learning}

\author{Yang Li\fnref{label1}}
\ead{liyangbx5433@163.com}
\address[label1]{School of Automation, Nanjing University of Science and Technology, 200 Xiaolingwei Street, Nanjing 210094, Jiangsu Province, China}

\author{Shenglan Yuan\corref{cor1}\fnref{label2,label3}}
\ead{shenglanyuan@gbu.edu.cn}
 \cortext[cor1]{Corresponding author}
\address[label2]{School of Sciences, Great Bay University, Songshan Lake International Innovation Entrepreneurship Community A5, Dongguan 523000, Guangdong Province, China}
\address[label3]{Great Bay Institute for Advanced Study, Songshan Lake International Innovation Entrepreneurship Community A5, Dongguan 523000, Guangdong Province, China}

\author{Shengyuan Xu\fnref{label1}}
\ead{syxu@njust.edu.cn}

\begin{abstract}
Stochastic vegetation-water dynamical systems play a pivotal role in ecological stability, biodiversity, water resource management, and adaptation to climate change. This research proposes a machine learning-based method for analyzing rare events in stochastic vegetation-water dynamical systems with multiplicative Gaussian noise. Utilizing the Freidlin-Wentzell large deviation theory, we derive the asymptotic expressions for the quasipotential and the mean first exit time. Based on the decomposition of vector field, we design a neural network architecture to compute the most probable transition paths and the mean first exit time for both non-characteristic and characteristic boundary scenarios. The results indicate that this method can effectively predict early warnings of vegetation degradation, providing new theoretical foundations and mathematical tools for ecological management and conservation. Moreover, the method offers new possibilities for exploring more complex and higher-dimensional stochastic dynamical systems.
\end{abstract}




\end{frontmatter}

\section{Introduction}\label{I}
The stochastic noise is indeed very common in ecosystems, which refers to random fluctuations in ecological systems that are not predictable or explainable by deterministic factors \cite{QY,YLZ,TTYBD}. These fluctuations can arise from a variety of sources, including environmental conditions, dispersal patterns and interactions among species, etc. Under random fluctuations, rare events \cite{SSCA}, i.e., transition from one stable state to another, frequently occur, even for weak noise. Thus they are well worth investigating, such as vegetation degradation and species extinction \cite{YW}.

The Freidlin-Wentzell large deviation theory \cite{FW} is actually a principle that is mainly applied to study rare events of dynamical systems  with small random perturbations. As the noise intensity decreases, the convergence rate at which the sample trajectories converge to the reference orbit is exponential in terms of noise. Large deviation techniques are widely used, which have become an extremely active branch of applied probability. They can estimate the escape probability of stochastic systems, reckon the probability of deviation from the reference orbit, and quantify the asymptotic probability of errors in hypothesis testing.

In large deviation theory, the action functional is an important concept that builds the relationship between the stationary probability distribution and the path distribution of stochastic dynamical systems \cite{FW}. It is defined as the integral of a certain Lagrangian, similar to classical mechanics. The minimization of the action functional leads to the most probable path connecting given initial and final states, along which the stochastic system moves with highest probability than other paths \cite{YD,LDLZ}. Although variables in statistical mechanics do not move along a definite trajectory, the most probable path offers us a very intuitive way to comprehend the stochastic behaviors of the dynamical systems and predict their evolution. Furthermore, the most probable path can also be used in optimization and control problems to find optimal solutions that maximize or minimize the objective functions in stochastic dynamical systems, which provides an appropriate and effective method for exploring the properties of stochastic models in practical applications \cite{LKBMM}. The calculation of the most probable path typically involves statistical methods and numerical simulation techniques. For instance, the Monte Carlo method is a random sampling technique, which can be utilized to simulate the behavior of a system under diverse conditions and parameters. Through extensive simulation and statistics, we can more accurately find the most probable path and its associated probability \cite{DMSSS}.

In addition, quasipotential is commonly used in physics and engineering to generalize the concept of potential function to nongradient dynamical systems \cite{FW,ZAAH,Ao}. It is defined as the global minimal value of action functional about both possible paths and time length, which can be used to describe rare transition events in various physical systems, such as fluid mechanics, electrodynamics and ecosystems. The quasipotential function exponentially dominates the magnitude of mean first exit time and stationary probability density.

Mean first exit time is a statistical quantity that represents the expected value of the time required for a stochastic system to escape from a confined state or region for the first time. In stochastic processes and diffusion theory, it provides a measure of how long a  stochastic system can transition from one state to another, which can be valuable in grasping and optimizing the performance of complex dynamical system related to the properties of the medium, such as its diffusivity, geometry, and boundary conditions. In practical applications, it is a critical component used to assess the engineering reliability of first passage failure \cite{CZ} and to describe the activation process in neuron systems \cite{FTPVB}. Apart from the quasipotential, more accurate perturbation expression of mean first exit time also depends on the exponential prefactor function of WKB approximation \cite{NKMS,MST,MS,LZXDL}.

The traditional numerical methods for calculating these previously mentioned quantities include the geometric minimum action method (GMAM) \cite{HVE} and ordered upwind method (OUM) \cite{C}. The former method is grounded in geometric action, aiming to iterate the most probable paths with fixed endpoints in the path space and obtain their corresponding action functionals. It converts the time parameterization of the action functional into arc length, simplifying the integration over infinite time to a finite-length integration, thereby finding the path corresponding to the minimum action. OUM is a numerical method for discretizing phase space by considering the influence of each node in an ordered manner and using an upwind difference scheme to compute the quasipotential at each node. This method facilitates efficient computation within the discretized phase space. GMAM and OUM are two different numerical methods based on the concepts of geometric action and discretization, respectively, for solving different types of problems. However, both GMAM and OUM have nonnegligible shortcomings, such as multiple minima for GMAM and great computational cost for OUM, especially for the high-dimensional case. Due to their limitations, there is a need to explore new methods.

Machine learning is an important direction in the field of artificial intelligence \cite{A}. It analyzes large amounts of data and improves existing algorithms to make them more intelligent and enhance their generalization capabilities \cite{SWS}. The powerful features of machine learning are mainly manifested in the following aspects. Firstly, machine learning is data-driven. By learning from a large amount of data, it can discover patterns and rules in the data, enabling more accurate predictions and decisions \cite{B}. Secondly, machine learning algorithms can automatically adjust parameters and models, reducing human intervention and errors; thus, making the results more objective and precise \cite{HLWFKB}. Thirdly, machine learning algorithms can quickly adapt to different datasets and tasks, exhibiting strong adaptability, facilitating learning and application across diverse environments \cite{AR}. Fourthly, as algorithms continuously improve and optimize, the interpretability of machine learning models is also steadily improving, enabling a better understanding of the internal mechanisms of data and models \cite{CCC}. Fifthly, machine learning has widespread applications in various fields, such as image recognition, speech recognition, natural language processing, recommendation systems, financial risk control, etc., bringing about tremendous changes and adding significant business value  across different sectors \cite{BB}. In summary, machine learning is a powerful technology that facilitates more accurate and effective data processing through automatic learning and optimization algorithms \cite{JM}. With the continuous development of technology, the application prospects of machine learning will become even broader.

At present, machine learning has been widely applied in the study of stochastic dynamics. For instance, researchers designed some data-driven methods to discover stochastic dynamical systems with Gaussian or non-Gaussian noise \cite{LD,CYDK}. Xu \emph{et al}. developed a novel deep learning method to compute the probability density by solving the Fokker-Planck equation \cite{XZLZLK}. There are also some scholars devoted to exploring data-driven machine learning methods for computing the most probable paths of stochastic dynamical systems \cite{LDL,WGCD}. Based on large deviation theory and machine learning, Li \emph{et al}. \cite{LYX} calculated most probable transition path and designed a control strategy to control the mean first exit time  of stochastic dynamical systems to achieve a desired value. Machine learning can also be used to compute the large deviation prefactors in the highly complex nonlinear stochastic systems\cite{LYLL}. These methods to study rare events are mainly used for stochastic systems with additive noise.

The investigation of stochastic vegetation-water ecosystems is of great significance due to their notable impact on ecological stability, biodiversity, water resource management, adaptation to climate change, and soil health. Considering their complex structures with multiplicative noise \cite{ZXLQ}, this paper aims to develop a machine learning method to handle rare transition events of stochastic dynamical systems with multiplicative Gaussian noise, rather than the case with additive noise \cite{LZXDL,LYLL,LLR}. Furthermore, this method will be utilized to calculate the most probable path, quasipotential, and mean first exit time of the stochastic vegetation-water system \cite{ZXLQ}. The analysis of the dynamic behavior of this system provides a theoretical basis and mathematical methods for understanding and controlling the phenomenon of vegetation degradation.

The structure of this work is as follows. In Section \ref{SVW}, we present the framework of the stochastic vegetation-water system and explain the meaning of the involved parameters. We reveal the dynamical structures of the corresponding deterministic system and simulate the transition phenomenon between the metastable states of stochastic system. In Section \ref{TM}, we describe the concepts of large deviation theory, quasipotential and mean first exit time, and derive their expressions asymptotically. Then a machine learning method is proposed to compute these quantities for the vegetation-water system in both non-characteristic and characteristic boundary cases. In Section \ref{R}, we use the machine learning algorithm to compute the rare transition events of the stochastic vegetation-water system and provide a mathematical basis for early warning of vegetation degradation. In Section \ref{C}, we summarize the results of the paper and discuss some important future prospects.

\section{Stochastic vegetation-water system}\label{SVW}
The vegetation-water dynamical system is a complex ecosystem where the vegetation and water interact to maintain the balance.
However, due to various influencing factors such as climate change, terrain, soil type, and the hydrological cycle, the dynamics of the vegetation-water system often exhibit uncertain and unpredictable stochasticity that manifests in several ways. Firstly, meteorological processes, such as rainfall intensity and evaporation, are stochastic, leading to stochasticity in hydrological processes (e.g., river discharge, groundwater levels). These hydrological processes directly influence  the state of the vegetation. Secondly, the growth and distribution of vegetation are also affected by many stochastic factors. For example, seed dispersal, plant growth rates, and the occurrence of pests and diseases are stochastic processes that lead to stochasticity in vegetation spatial distribution and density. In addition, the interaction between vegetation and water in the system is also stochastic. For instance, vegetation consumes water through transpiration, and water returns to the system through soil infiltration, surface runoff, etc.
To understand and predict the stochasticity of the vegetation-water dynamical system, mathematical tools such as probability theory and stochastic processes are required to establish quantitative models of the system's state and behavior. Furthermore, modern technological means, like machine leaning techniques, are employed to obtain extensive observational data of the system, supporting model validation and refinement. The stochasticity of vegetation-water systems is a crucial characteristic that plays a key role in predicting and responding to practical issues such as water resource management and ecological restoration. Studying and comprehending this stochasticity is essential for effective management and conservation efforts.

More specifically, we consider a stochastic vegetation-water system
\begin{equation}\label{vw}
\begin{array}{l}
\dot{x}_1=\rho x_1\big(x_2-\frac{x_1}{K}\big)-\beta\frac{x_1}{x_1+x_0}+\sigma_{1}x_1^{2}\xi_{1}(t),\\
\dot{x}_2=R-\alpha x_2-\lambda x_1x_2+\sigma_{2}\xi_{2}(t),
\end{array}
\end{equation}
where the variable $x_1$ represents the biomass of vegetation, while $x_2$ signifies the moisture level of the soil.
 The interaction terms $\rho x_1x_2$ and $-\lambda x_1x_2$ depict the relationship between vegetation biomass and water.
The term $-\rho x_1^2/K$ limits the growth of biomass due to competition for shared resources, such as water or soil nutrients.
The term  $-\beta x_1/(x_1+ x_0)$ illustrates the impact of herbivores and other influencing factors. The parameter
$R$ denotes the average rainfall, and the term $-\alpha x_2$ stands for the loss of water from the soil, which could be caused by percolation or evaporation. Taking into account random environmental disturbances affecting vegetation competition for shared resources and rainfall, the term $-\rho x_1^2/K$  transforms into
$-\rho x_1^2(1+\tilde{\sigma}_{1}\xi_1(t))/K$ and $R$ becomes $R(1+\tilde{\sigma}_{2}\xi_2(t))$. Consequently, these disturbances are captured by multiplicative noise
$\sigma_{1}x_1^2\xi_1(t)$ and additive noise $\sigma_{2}\xi_2(t)$, where $\sigma_{1}=-\rho\tilde{\sigma}_{1}/K$ and $\sigma_{2}=R\tilde{\sigma}_{2}$. The driving terms $\xi_{1}(t)$ and $\xi_{2}(t)$ are independent Gaussian white noises with
\begin{equation*}
\mathbb{E}[\xi_{i}(t)]=0,\quad \mathbb{E}[\xi_{i}(t)\xi_{i}(t+\tau)]=\varepsilon \delta(\tau),\quad i=1,2.
\end{equation*}
Here, $\varepsilon$ is a small parameter, implying that the noise intensity is weak.
The drift coefficient and diffusion matrix can be written as
\begin{equation*}
b(x)=\left(
       \begin{array}{c}
         \rho x_1\big(x_2-\frac{x_1}{K}\big)-\beta\frac{x_1}{x_1+x_0} \\
          R-\alpha x_2-\lambda x_1x_2\\
       \end{array}
     \right),\quad a(x)=\sigma(x)\sigma^{\top}(x)=\left(
                          \begin{array}{cc}
                            \sigma_{1}^{2}x_1^{4} & 0 \\
                            0 &  \sigma_{2}^{2}\\
                          \end{array}
                        \right),
\end{equation*}
where
\begin{equation*}
 \sigma(x)=\left(
 \begin{array}{cc}
 \sigma_{1}x_1^{2} & 0 \\
 0 & \sigma_{2} \\
 \end{array}
 \right).
\end{equation*}
We first consider the corresponding deterministic system with $\sigma_{1}=\sigma_{2}=0$:
\begin{equation}\label{dvw}
\begin{array}{l}
\dot{x}_1=\rho x_1\big(x_2-\frac{x_1}{K}\big)-\beta\frac{x_1}{x_1+x_0},\\
\dot{x}_2=R-\alpha x_2-\lambda x_1x_2.
\end{array}
\end{equation}
In this paper, we fix these system parameters as $\rho=1$, $K=10$, $\beta=3$, $x_0=1$, $\alpha=1$, $\lambda=0.12$.

\begin{figure}
	\centering
	\includegraphics[width=7cm]{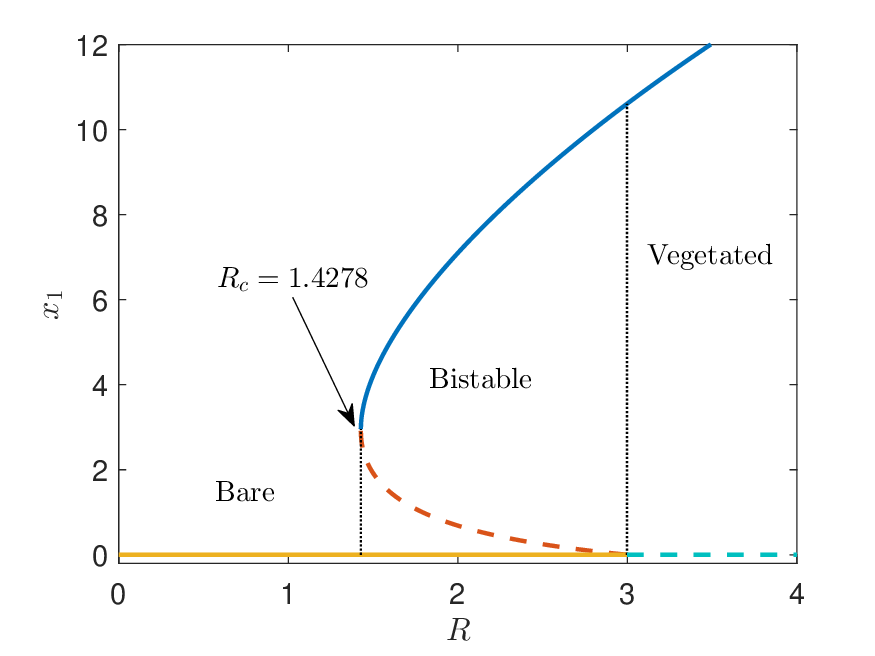}
	\caption{Bifurcation digram of the vegetation system about the parameter $R$. It exhibits two branches, the blue curve representing the stable equilibria, while the red dashed  branch representing the unstable equilibria. Upon varying the control parameter $R$, the two branches approach each other until they meet at the critical value of $R_c=1.4278$. At this point, they annihilate each other in a saddle-node bifurcation.}
	\label{fig1}
\end{figure}

Let $b(x)=\textbf{0}$ and we can derive the fixed points of the system \eqref{dvw}. It is observed that $\text{SN1}=\big(0,\frac{R}{\alpha}\big)$
is a trivial fixed point for arbitrary $R$. When $x_1\neq0$, $b(x)=\textbf{0}$ implies that
\begin{equation*}
\begin{array}{l}
f(x_1)=\rho\lambda x_1^{3}+\rho(\lambda x_0+\alpha)x_1^{2}+(\rho\alpha x_0+K\beta\lambda-\rho RK)x_1+K\beta\alpha-\rho RKx_0=0,\\
x_2=\frac{R}{\lambda x_1+\alpha}.
\end{array}
\end{equation*}
As shown in Fig. \ref{fig1}, a saddle-node bifurcation occurs at $R=R_c$. Therefore, a node SN2 and a saddle US emerge.
Let $f(x_1)=0$ and $f'(x_1)=0$. We have $R_c=1.4278$. SN1 is stable if $R<2.998$, otherwise it is unstable due to the collision of US.
Above all, there are three states for deterministic vegetation-water system  \eqref{dvw}  depending on  the value range of $R$:
bare for $R<R_c$, bistable for $R_c<R<2.998$, and vegetated for $R>2.998$.

It can be seen from the saddle-node bifurcation diagram in Fig. \ref{fig1} that the biomass of vegetation depends greatly on the average rainfall $R$. This finding aligns with the actual ecological phenomena. When $R<R_c$, the rainfall is insufficient for vegetation to survive, which corresponds to the bare state; when $R>2.998$, the rainfall is ample, and thus the vegetation can grow fully without the phenomenon of vegetation disappearance, corresponding to the vegetated state. Meanwhile, for rainfall values within the range $R_c<R<2.998$, the bistable phenomenon emerges as the rainfall is neither too little nor too much. Investigating the vegetation state within this range of rainfall has practical significance for applications.

\begin{figure}
	\centering
	\includegraphics[width=7cm]{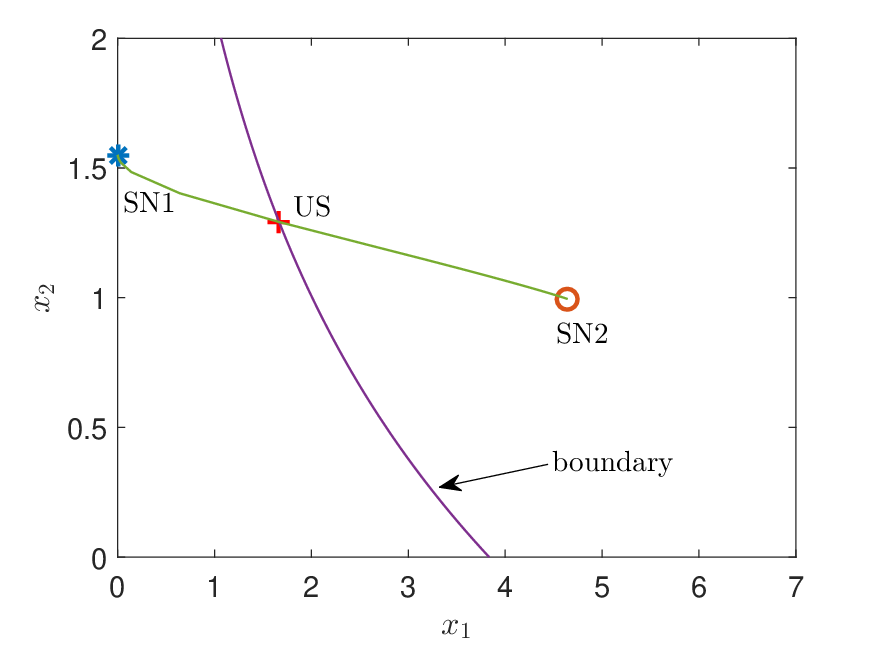}
	\caption{For the selected value of $R=1.55$, the system \eqref{dvw} possesses two stable fixed points $SN1=(0,1.55)$ and $SN2=(4.6366,0.9959)$, which are separated by the stable manifold of the saddle point US=(1.6667,1.2917). This stable manifold is delineated by a purple curve. Additionally, the unstable manifold of the saddle point US is denoted by a green curve.}
	\label{fig2}
\end{figure}

In this paper, we choose $R=1.55$, i.e., bistable state, for investigation. As demonstrated in Fig. \ref{fig2}, the system \eqref{dvw} has two stable fixed points
$\text{SN1}=(0,1.55)$ and $\text{SN2}=(4.6366,0.9959)$. The basins of attraction for these fixed points are separated by the stable manifold of the saddle $\text{US}=(1.6667,1.2917)$, denoted
by a purple curve. The unstable manifold of US is indicated by a green curve.

\begin{figure}
	\centering
	\includegraphics[width=7cm]{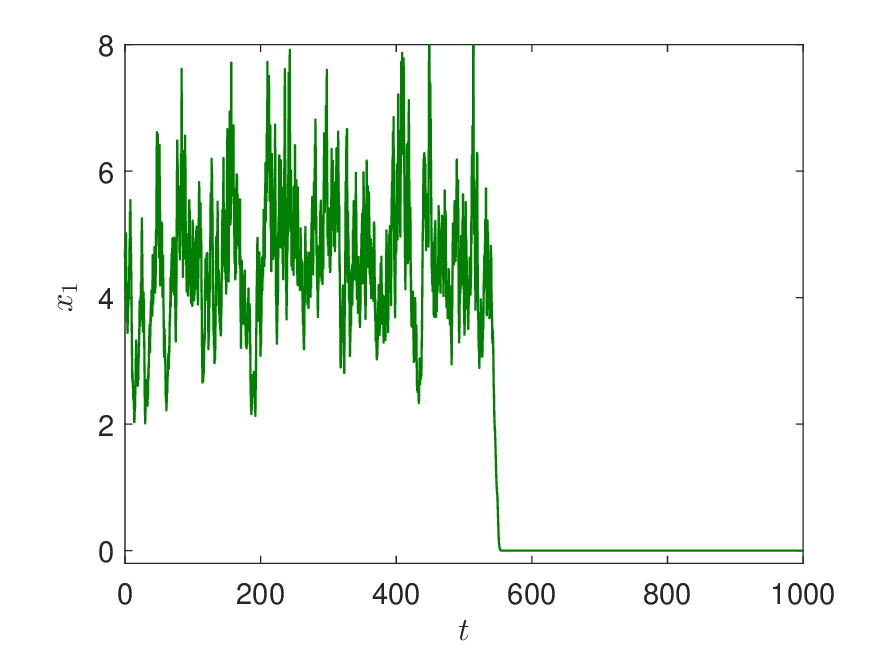}
	\caption{A representative transition trajectory of the stochastic vegetation-water model \eqref{vw} depicted through Monte Carlo simulation.}
	\label{fig3}
\end{figure}

Based on the ecological significance of stochastic vegetation-water model \eqref{vw}, we investigate the transition phenomena of the system initially located at SN2, i.e., with lush vegetation, under random perturbations approaching SN1. Since the noise intensity $\varepsilon$ is small, these phenomena are referred to as rare events. Indeed, when the system crosses the boundary, i.e., the stable manifold of US, it will flow along the unstable manifold of the US to the point SN1. Therefore, we mainly focus on the process of the system escaping from the attractor domain of SN2 driven by noise perturbations, which can provide important information for early warning and intervention of vegetation loss.
A typical noise induced transition trajectory simulated by Monte Carlo is illustrated in Fig. \ref{fig3}.

\section{Theory and method}\label{TM}
Under the condition of weak Gaussian noise, rare exit events can be analyzed by utilizing Freidlin-Wentzell large deviation theory.  In this section, we first review the concepts of the most probable exit path, quasipotential, and mean first exit time, and derive their asymptotic expressions. Then we design a machine learning method to compute these quantities numerically.

\subsection{Large deviation theory}\label{LDT}
Now we reformulate the stochastic vegetation-water system \eqref{vw} into the following form
\begin{equation*}
dx(t)=b(x)dt+\sigma(x)dB^{\varepsilon}(t),
\end{equation*}
where $B^{\varepsilon}(t)=(B_1^{\varepsilon}(t),B_2^{\varepsilon}(t))^{\top}$ is a two-dimensional Brownian motion. Given the scenario of weak noise intensity,
the stationary distribution $p_s(x)$ of the stochastic system can be assumed to have the following WKB form
\begin{equation*}
p_s(x)\sim C(x)\exp\{-\varepsilon^{-1}V(x)\},
\end{equation*}
where $C(x)$ is referred to as exponential prefactor, and $V(x)$ is called quasipotential, which characterizes the possibility of the stochastic state fluctuating into
the vicinity of the specific point $x$. In addition to the stationary distribution, the quasipotential also exponentially dominates the magnitude of mean first exit time and exit location distribution. According to Freidlin-Wentzell large deviation theory, the quasipotential is defined as the minimum of the action functional along the absolutely continuous path
connecting the fixed point $\bar{x}$ and the specific point $x$, in the sense that
\begin{equation*}
V(x):=\inf_{T>0}\inf_{\varphi\in C[0,T]}\{S(\varphi):\varphi(0)=\bar{x},\, \varphi(T)=x\},
\end{equation*}
where the action functional $S(\varphi)$ has the following form
\begin{equation*}
S(\varphi)=\frac{1}{2}\int_0^T(\dot{\varphi}-b(\varphi))^{\top}a^{-1}(\dot{\varphi}-b(\varphi))dt.
\end{equation*}

By substituting the WKB approximation into the stationary Fokker-Planck equation and combining the lowest order terms of $\varepsilon$, we obtain Hamilton-Jacobi equation
\begin{equation*}
\langle \nabla V(x),b(x)\rangle+\frac{1}{2}\langle \nabla V(x),a(x)\nabla V(x)\rangle=0.
\end{equation*}
Note that the above equation has a geometric interpretation. By combining the two terms on the left, we gain
\begin{equation*}
\langle \nabla V(x),b(x)+\frac{1}{2}a(x)\nabla V(x)\rangle=0.
\end{equation*}
Therefore, we get the following orthogonal relationship
\begin{equation*}
 \nabla V(x)\perp b(x)+\frac{1}{2}a(x)\nabla V(x).
\end{equation*}
Define $l(x):=b(x)+\frac{1}{2}a(x)\nabla V(x)$. The vector field $b(x)$ has the following  decomposition
\begin{equation*}
b(x)=-\frac{1}{2}a(x)\nabla V(x)+l(x).
\end{equation*}

The next-to-leading-order approximation of the WKB approximation to the stationary Fokker-Planck equation leads to
the following transport equation about the exponential prefactor $C(x)$:
\begin{equation*}
\langle\nabla C, b+a\nabla V\rangle+C\big(\text{div}b+\frac{1}{2}a:H(V)+\langle A,\nabla V\rangle\big)=0,
\end{equation*}
where $H$ denotes Hessian matrix, $a:H=\sum_{i,j}a_{ij}H_{ij}$, and the vector $A(x)$ is defined by $A_i(x)=\sum\limits_{j=1}^{2}\frac{\partial a_{ij}}{\partial x_j}$, i.e.,
\begin{equation*}
A(x)=\left(
      \begin{array}{c}
        4\sigma_1^2x_1^3 \\
        0\\
      \end{array}
    \right).
\end{equation*}
Putting the decomposition $b(x)=-\frac{1}{2}a(x)\nabla V(x)+l(x)$ into the transport equation yields the following result
\begin{equation*}
\langle\nabla C,\frac{1}{2}a\nabla V+l\rangle+C(\text{div}l+\frac{1}{2}\langle A,\nabla V\rangle)=0.
\end{equation*}
The above equation can be rewritten as
\begin{equation}\label{te}
\langle\nabla\ln C,\frac{1}{2}a\nabla V+l\rangle=-F,
\end{equation}
where
\begin{equation*}
F(x)=\text{div}l(x)+\frac{1}{2}\langle A(x),\nabla V(x)\rangle.
\end{equation*}
According to Freidlin-Wentzell large deviation theory,
the most probable path of fluctuation dynamics satisfies
the following equation
\begin{equation*}
\dot{\varphi}^x(t)=b(\varphi^{x}(t))+a(\varphi^{x}(t))\nabla V(\varphi^{x}(t))=\frac{1}{2}a(\varphi^{x}(t))\nabla V(\varphi^{x}(t))+l(\varphi^{x}(t)).
\end{equation*}
This can also be confirmed by the results of the method of characteristics applied to the Hamilton-Jacobi equation. Therefore,
along the most probable path, the left-hand side of equation \eqref{te} is transformed into the complete differential of $\ln C(x)$, i.e.,
\begin{equation*}
\frac{d}{dt}\ln C(\varphi^x(t))=-F(\varphi^x(t)).
\end{equation*}
So that we can integrate the prefactor function as
\begin{equation}\label{Cx}
C(x)\sim\exp\{-\int_{-\infty}^{0}F(\varphi^x(t))dt\}.
\end{equation}
If $x$ is a saddle point, then the upper limit of the above integral is infinite.

Assume that $D$ is the attraction domain of SN2. Define the first exit time of the stochastic state from $D$ as
\begin{equation*}
\tau^{\varepsilon}_{D}=\inf\{t\geq0:x(t)\notin D, x(0)=\text{SN2}\}.
\end{equation*}
It is a random time, and its average quantity is called the mean first exit time. This can provide important quantitative information  regarding the disappearance of vegetation.
According to large deviation theory, the mean first exit time is exponentially dominated by the minimal value of the quasipotential along the boundary $\partial D$:
\begin{equation*}
\lim_{\varepsilon\rightarrow0}\varepsilon\ln\mathbb{E}\tau^{\varepsilon}_{D}=\inf_{x\in\partial D}V(x).
\end{equation*}
Usually, this minimization is achieved at the saddle US. Then we have
\begin{equation*}
\mathbb{E}\tau^{\varepsilon}_{D}=L^{\varepsilon}_{D}\exp\{V(\text{US})\}.
\end{equation*}

In general, We can calculate the prefactor $L_{D}^{\varepsilon}$ in two distinct scenarios.

\textbf{Case A.} non-characteristic boundary (see \cite{BR22,BR} for reference).

We assume that the domain $D$ is an open, smooth and connected subset of $\mathbb{R}^{n}$ that satisfies the following conditions, where $n(y)$ represents the exterior normal vector at $y\in\partial D$.
\begin{description}
  \item[(A1)] The deterministic system $\dot{x}=b(x)$ has a unique fixed point $\bar{x}$ within $D$ that attracts all the trajectories originating from $D$. Additionally, the inner product between $b(y)$ and $n(y)$ is negative for all $y$ on the boundary $\partial D$, i.e., $\langle b(y),n(y)\rangle<0$ for all $y\in \partial D$.
  \item[(A2)] The function $V$ is continuously differentiable ($C^1$) in $D$; for any $x\in\bar{D}$, the most probable path $\varphi_{t}^{x}$ approaches $\bar{x}$ as $t\rightarrow-\infty$; and $\langle \frac{1}{2}a(y)\nabla V(y)+l(y),n(y)\rangle >0$ for all $y\in \partial D$.
  \item[(A3)]  The minimum of $V$ over $\partial D$ is attained at a single point $x^{*}$. At this point,
  \begin{equation*}
  \mu^{*}=\langle \frac{1}{2}a(x^{*})\nabla V(x^{*})+l(x^{*}),n(x^{*})\rangle>0,
  \end{equation*}
  and the quadratic form $h^{*}: \xi\mapsto\langle\xi,\nabla^{2}V(x^{*})\xi\rangle$ has positive eigenvalues on the hyperplane $n(x^{*})^{\perp}=\{\xi\in\mathbb{R}^{n}: \langle\xi,n(x^{*})\rangle=0\}$.
\end{description}

If $\langle b(y),n(y)\rangle<0$ for all $y\in \partial D$, the boundary is designated as non-characteristic. This condition guarantees
that the dynamical trajectories, originating from the closure $\bar{D}$, will remain confined within $D$, and that the vector field is perpendicular to the boundary. Leveraging Assumption (A1), we derive the integral formula
\begin{equation*}
\lambda_{D}^{\varepsilon}=\int\limits_{x\in \partial D}\langle \frac{1}{2}a(x)\nabla V(x)+l(x),n(x)\rangle C(x)\exp\{-\varepsilon^{-1}V(x)\}dx
\end{equation*}
for the exit rate $\lambda_{D}^{\varepsilon} = [\mathbb{E} \tau_{D}^{\varepsilon}] ^{-1}$. By employing the second-order expansion of the potential $V$ in the vicinity of the point $x^{*}$, we obtain an equivalent relation for the prefactor
\begin{equation}\label{LcaseA}
\begin{aligned}
L_{D}^{\varepsilon}&\sim\frac{1}{C(x^{\ast})\mu^{\ast}}\sqrt{\frac{\det h^{\ast}}{(2\pi\varepsilon)^{n-1}}}  \\
                   &\sim\frac{1}{\mu^{\ast}}\sqrt{\frac{\det h^{\ast}}{(2\pi\varepsilon)^{n-1}}} \exp\{\int_{\infty}^{0}F(\varphi^{x^{\ast}}(t))dt\},
\end{aligned}
\end{equation}
where we utilize the approximation in \eqref{Cx}.

\textbf{Case B.} characteristic boundary \cite{BR22,BR}.

The basin domain $D$ is characteristic in the sense that the inner product between the vector field $b(y)$ and the exterior normal vector $n(y)$ vanishes for all points $y$ within the domain $D$, i.e., $\langle b(y),n(y)\rangle=0$ for all $y\in D$. We consider the metastable scenario that the deterministic system $\dot{x}=b(x)$ possesses two stable fixed points $\bar{x}_1$ and $\bar{x}_2$, whose respective basins of attraction are divided by a smooth hypersurface $S$. We focus on exit events from the basin of attraction $D$ associated with $\bar{x}_1$ and introduce the following set of assumptions.
\begin{description}
  \item[(B1)] All trajectories of the deterministic system $\dot{x}=b(x)$ initiated on the hypersurface $S$ remain confined to $S$ and
  ultimately converge to a single fixed point $x^{\ast}\in S$; furthermore, the Jacobi matrix $\nabla b(x^{\ast})$ possesses $n-1$ eigenvalues with negative real part and a single positive eigenvalue denoted by $\lambda^{*}$.
  \item[(B2)] With respect to the quasipotential $V$ associated with $\bar{x}_1$, there exists a unique (up to time shift) trajectory $\rho=(\rho_t)_{t\in\mathbb{R}}\subset D$ such that
  \begin{equation*}
  \lim_{t\rightarrow-\infty}\rho_t=\bar{x}_{1},\quad \lim_{t\rightarrow+\infty}\rho_t=x^{*},\quad\text{and}\quad V(x^{*})=\mathcal{S}_{-\infty,+\infty}[\rho].
  \end{equation*}
  \item[(B3)] The quasipotential $V$ is smooth in a neighborhood of $\rho=(\rho_t)_{t\in\mathbb{R}}$. Additionally, the vector field $l$ defined by  $l(x)=b(x)+\frac{1}{2}a(x)\nabla V(x)$ satisfies the orthogonality relation $\langle \nabla V(x),l(x)\rangle=0$.
\end{description}

In this context, the quasipotential $V$ attains its minimum on the hypersurface $S$ precisely at the point $x^{*}$. Moreover, the trajectory $\rho$ is designated as the most probable exit path, satisfying the differential equation
\begin{equation*}
\dot{\rho}_t=\frac{1}{2}a(\rho_t)\nabla V(\rho_t)+l(\rho_t),\quad\forall t\in\mathbb{R}.
\end{equation*}
For any given $t\in\mathbb{R}$, this path coincides with the trajectory $(\varphi_s^{x})_{s\leq0}$ that connects $\bar{x}$ to $x=\rho_t$, according to the relation
\begin{equation*}
\varphi_{s}^{x}=\rho_{s+t},\quad\forall s\leq0.
\end{equation*}
To describe the prefactor $L_D^{\epsilon}$ in this scenario, we formulate the following supplementary assumption.
\begin{description}
  \item[(B4)] The matrix $H^{\ast}=\lim_{t\rightarrow+\infty}\nabla^{2}V(\rho_t)$ exists and possesses $n-1$ positive eigenvalues and a single negative eigenvalue.
\end{description}

Relying on these four assumptions, an asymptotic formula for estimating the expected time taken by the process to exit the domain $D$ is given by
\begin{equation}\label{LcaseB}
L_{D}^{\varepsilon}\sim\frac{\pi}{\lambda^{\ast}}\sqrt{\frac{|\det H^{\ast}|}{\det\nabla^{2}V(\bar{x})}}\exp\Big\{\int_{-\infty}^{\infty}F(\varphi(t))dt\Big\}.
\end{equation}
This expression provides an approximation for the expected exit time from the domain $D$
under the given conditions.

\subsection{Machine learning algorithm}\label{MLA}
It is seen in subsection \ref{LDT} that the computations of the most probable path and the mean first exit time depend on
the results of quasipotential and rotational component, i.e., the decomposition of the vector field. In this subsection,
we aim to propose a machine learning method to compute these quantities of the stochastic vegetation-water system \eqref{vw} based on
this decomposition.

We design a neural network architecture to achieve this goal. The input of the network is set as the
coordinate $x=(x_1,x_2)^{\top}$. The output of the network is $(\hat{V}_\theta,l_\theta)\in\mathbb{R}^{n+1}$.
The quasipotential function is defined as $V_\theta(x)=\hat{V}_\theta(x)+|x-\bar{x}|^{2}$ to guarantee
the unboundedness of $V(x)$ and $|\nabla V(x)|$. Here, $\bar{x}$ is the stable fixed point SN2, and $\theta$
denotes the training parameters of the neural network.

In order to train the neural network, we choose $N$ points randomly in the attraction domain of SN2 and construct
a loss function as follows:
\begin{equation*}
L=L_{\text{dyn}}+\lambda_{1}L_{\text{orth}}+\lambda_{2}L_{0}.
\end{equation*}
Since the vector field has the decomposition $b(x)=-\frac{1}{2}a(x)\nabla V(x)+l(x)$,
the first part $L_{\text{dyn}}$ of the loss function can be set as
\begin{equation*}
L_{\text{dyn}}=\frac{1}{N}\sum_{i=1}^{N}[b(x_i)+\frac{1}{2}a(x_i)\nabla V_\theta(x_i)-l_{\theta}(x_i)]^2,
\end{equation*}
where the gradient of the quasipotential is realized by automatic differentiation technique. Due to
the orthogonal relation $\nabla V(x)\bot l(x)$, the second part $L_{\text{orth}}$ of the loss function
is assigned as
\begin{equation*}
L_{\text{orth}}=\frac{1}{N}\sum_{i=1}^{N}\frac{[\nabla V_{\theta}(x_i)\cdot l_{\theta}(x_i)]^2}{|\nabla V_{\theta}(x_i)|^2|l_{\theta}(x_i)|^2+\delta}.
\end{equation*}
Here, the small parameter $\delta\ll1$ is chosen to avoid a zero denominator, ensuring numerical stability. Besides, the third part $L_0$ of the loss function is set as
$L_0=V_\theta(\bar{x})^{2}$ to guarantee the fact that the quasipotential of the stable fixed point $\bar{x}=\text{SN2}$ is zero.
In this paper, we choose the weight parameters $\lambda_{1}=1$ and $\lambda_{2}=0.1$, which are used to balance the three parts of loss function.

After training the neural network, we obtain the quasipotential function $V_{\theta}(x)$ and the rotational component $l(x)$. Thus the most probable path can be integrated by the equation
\begin{equation*}
\dot{x}=b(x)+a(x)\nabla V_{\theta}(x)
\end{equation*}
in reverse time, starting from the end point. Additionally, the mean first exit time can be computed using the asymptotic expression provided in subsection \ref{LDT}.

\section{Results}\label{R}
In this section, we present the results of applying the proposed machine learning algorithm to the stochastic vegetation-water system \eqref{vw}. The hyperparameters we set are as follows:
The neural network is configured with 6 hidden layers, each having 20 nodes. The activation function in hidden layers is chosen as $\tanh$, while the output layer uses the identify function. We
utilize the Adam optimizer with a learning rate of 0.001. The small parameter in the loss function is set to $\delta=0.001$.
The neural network is trained for 1000000 epochs.

In Section \ref{SVW}, upon observing the expressions of the noise intensities $\sigma_{1}=-\rho\tilde{\sigma}_{1}/K$ and $\sigma_{2}=R\tilde{\sigma}_{2}$ with $\rho=1$, $K=10$, and $R=1.55$, it is worth noting that the order of magnitude of $\sigma_1$ should be significantly smaller than that of $\sigma_2$ in practical system. Consequently, for this context, we assume the values of $\sigma_1=0.1$ and $\sigma_2=1$.

\begin{figure}
	\centering
	\includegraphics[width=7cm]{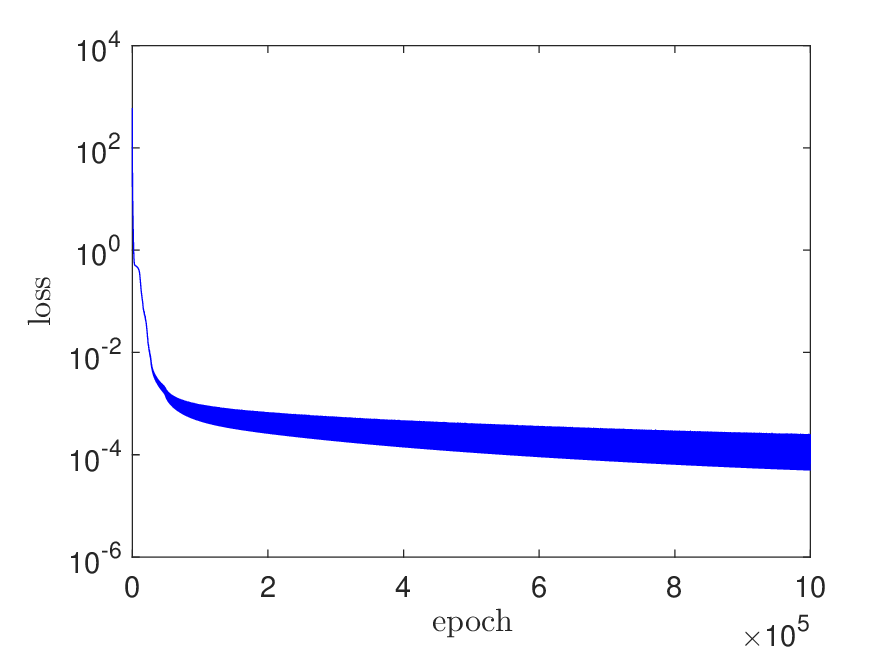}
	\caption{The loss function is decreasing with respect to the increasing number of epochs during the training of the neural network.}
	\label{fig4}
\end{figure}

\begin{figure}
	\centering
	\subfigure[Learned $V_{\theta}(x)$]
	{\includegraphics[width=4.4cm]{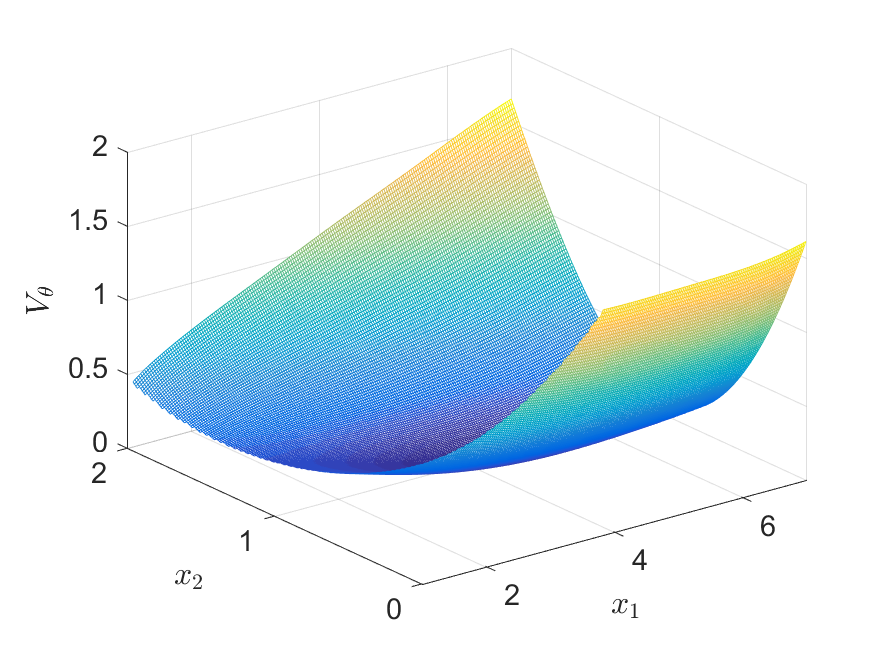}}
	\subfigure[Learned $l_{1\theta}(x)$]
	{\includegraphics[width=4.4cm]{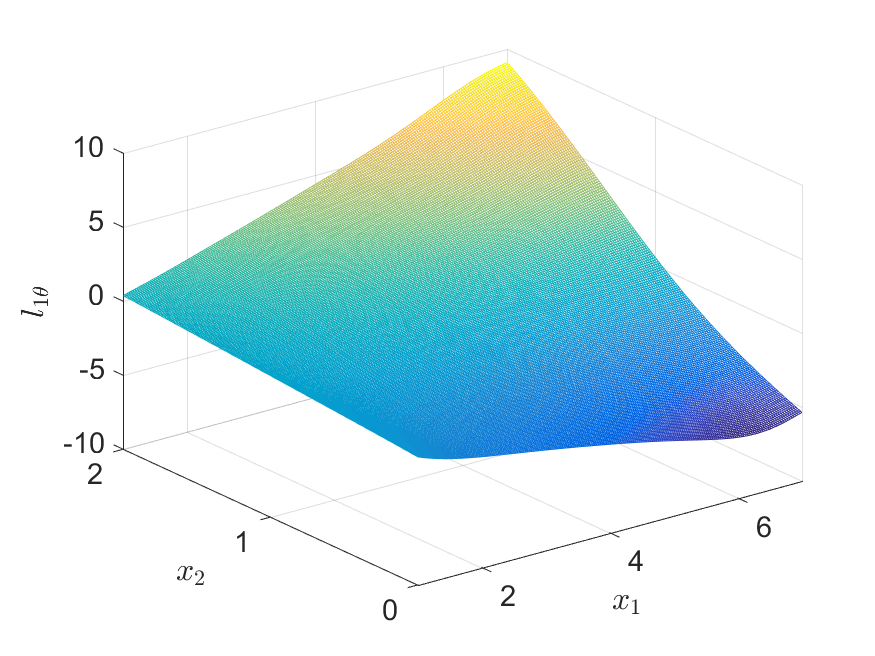}}
	\subfigure[Learned $l_{2\theta}(x)$]
	{\includegraphics[width=4.4cm]{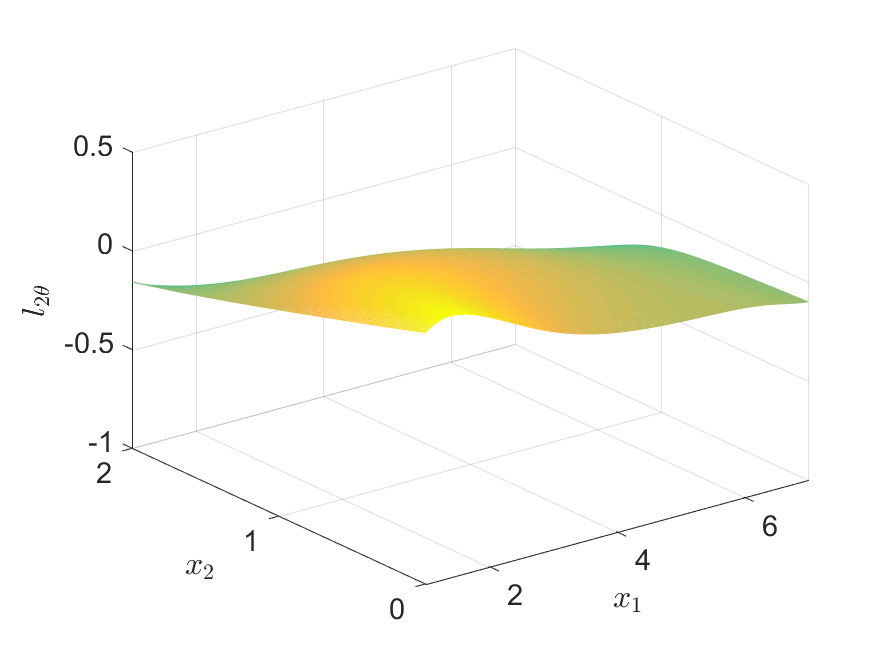}}
	\caption{The quasipotential function $V_\theta(x)$
and the two rotational components $l_{1\theta}(x)$, $l_{2\theta}(x)$ are learned by machine learning.}
	\label{fig5}
\end{figure}

Now we apply our proposed method to the exit problem of the system \eqref{vw}. We randomly and uniformly select 10000 points in the domain $[1,7]\times[0,2]$, within which $N=8018$ collocation points located on the right-hand side of the stable manifold of US are used to train the neural network. As depicted in Fig. \ref{fig4}, the loss function is reduced to the magnitude of $10^{-5}$, indicating good convergence of the algorithm.  The quasipotential function and rotational components of the training results are exhibited in Fig. \ref{fig5}.

We first consider the case of a non-characteristic boundary. We take $x_1=3$ as the boundary, as it has practical ecological significance. On the one hand, this boundary is located on the right side of the natural boundary, i.e., the stable manifold of the saddle point, and thus it can be used for early warning of vegetation degradation, leaving sufficient time for manual intervention. On the other hand, this boundary is relatively simple, only requiring measurement of vegetation biomass.

\begin{figure}
	\centering
	\includegraphics[width=7cm]{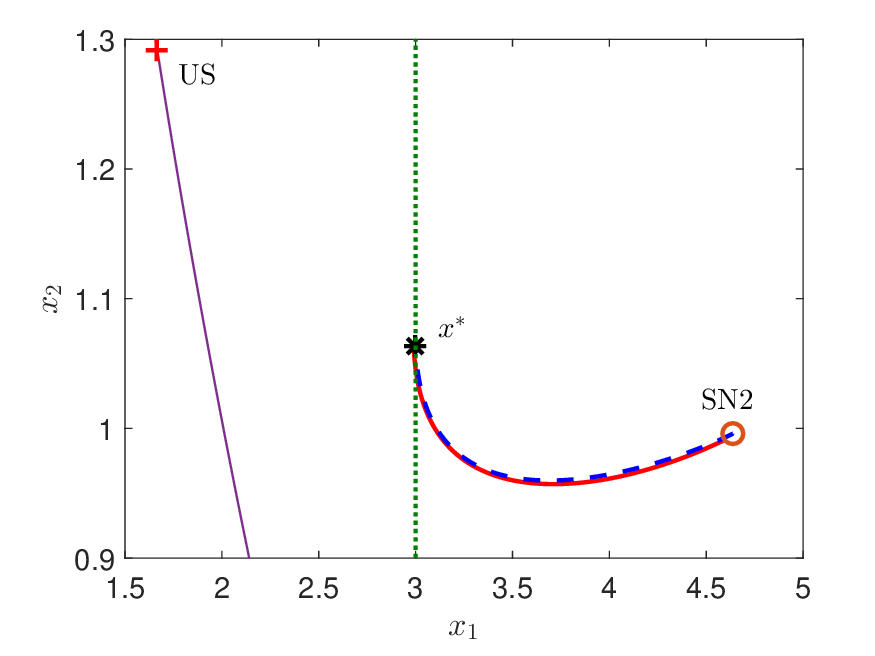}
	\caption{Along with the identification of the point $x^{\ast}=(3,1.0632)$ as the location of minimal quasipotential, the results of the most probable path connecting $x^{\ast}$ and SN2 are consistent with the outcomes of the machine learning.}
	\label{fig6}
\end{figure}

By employing the gradient descent method, we can locate the point $x^{\ast}=(3,1.0632)$ with the minimal quasipotential
on the boundary $x_1=3$. Starting from $x^{\ast}$, we can use the inverse time integration to determine the most probable path connecting $x^{\ast}$ and SN2, as demonstrated by the red curve in Fig. \ref{fig6}. The blue dashed line represents the path obtained by the shooting method \cite{BMLSM}. The consistency of these results with those obtained through machine learning is evident.

\begin{figure}
	\centering
	\includegraphics[width=7cm]{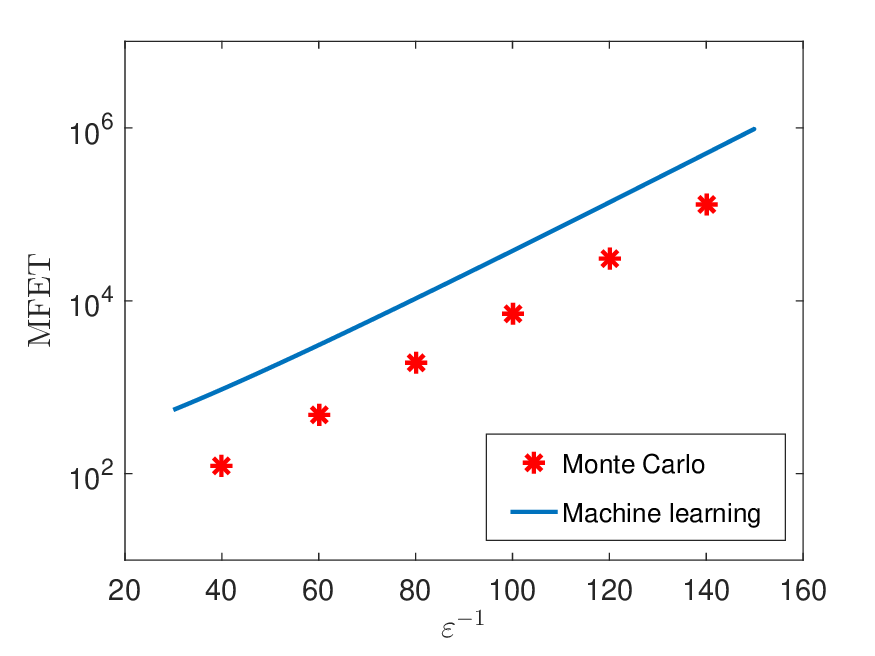}
	\caption{The machine learning result, represented by the blue curve, is in agreement with the Monte Carlo simulation, indicated by the red stars. This agreement  underscores the validity of both methods and their consistency in predicting the relationship between mean first exit time and noise intensity, despite the presence of a small error in the non-characteristic boundary case.}
	\label{fig7}
\end{figure}

Denote $\bar{H}=\nabla^{2}V_{\theta}(\text{SN2})$.
We expand the Hamilton-Jacobi equation near the fixed point $\bar{x}=\text{SN2}$ to obtain the algebraic Riccati equation
\begin{equation*}
\bar{H}^{-1}\bar{Q}^{\top}+\bar{Q}\bar{H}^{-1}=\bar{a},
\end{equation*}
where
\begin{equation*}
\bar{Q}=-\nabla b(\bar{x}),\quad \bar{a}=a(\bar{x}).
\end{equation*}
Therefore, we can get
\begin{equation*}
\bar{H}=\left(
          \begin{array}{cc}
            0.0543 & 0.0608 \\
            0.0608 & 2.9133 \\
          \end{array}
        \right).
\end{equation*}
Then $\det(\bar{H})=0.1546$. Besides,
\begin{equation*}
\begin{array}{l}
\mu^{\ast}=-\big(\frac{1}{2}\sigma_{1}^{2}(x_{1}^{\ast})^{4}\frac{\partial V_\theta}{\partial x_1}(x^{\ast})+l_{1\theta}(x^{\ast})\big)=0.022, \\
\det(h^{\ast})=\frac{\partial^{2}V_\theta}{\partial x_2^{2}}(x^{\ast})=2.2602,\quad V_\theta(x^{\ast})=0.0691.
\end{array}
\end{equation*}
Due to the asymptotic expression of mean first exit time, we can compute the functional relation between mean first exit time and the noise intensity $\varepsilon$, as illustrated in Fig. \ref{fig7}. The Monte Carlo simulations confirm the machine learning results, wherein the small error mainly stems from the asymptotic expression itself.

Next, we consider the characteristic boundary, namely, the stable manifold of US. This serves as the natural boundary of the system \eqref{vw}. When noise perturbation causes the system to reach this boundary, it tends to move towards SN1 along the unstable manifold of US. Therefore, once the system escapes from this boundary, it is not far from vegetation degradation, and effective measures need to be taken to change the situation.

\begin{figure}
	\centering
	\includegraphics[width=7cm]{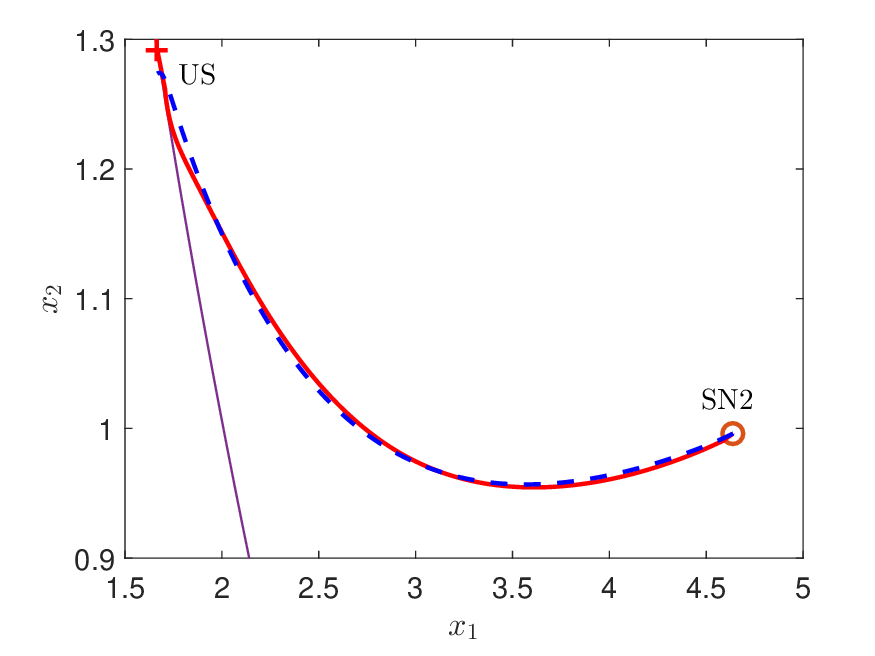}
	\caption{The machine learning-computed most probable path from SN2 to US closely agrees with the result obtained from the shooting method.}
	\label{fig8}
\end{figure}

According to large deviation theory, the saddle US represents the point with minimal quasipotential on the boundary i.e., the exit point. To obtain the most probable path, we integrate the equation $\dot{x}=b(x)+a(x)\nabla V_\theta(x)$, commencing from the neighborhood of US with time reversal. As seen in Fig. \ref{fig8}, the most probable path computed by machine learning aligns well with the result obtained via the shooting method. It should be noted that on the most probable path from SN2 to US, the soil moisture level $x_2$ undergoes a process of first decreasing and then increasing. Although the water content at US is much higher than that of the steady state SN2, vegetation inevitably degrades. Hence, in scenarios where vegetation is inadequate, enhancing soil moisture alone is insufficient to halt desertification. Under such circumstances, the effect of planting plants artificially is significantly more beneficial, which is far greater than merely increasing soil moisture.

\begin{figure}
	\centering
	\includegraphics[width=7cm]{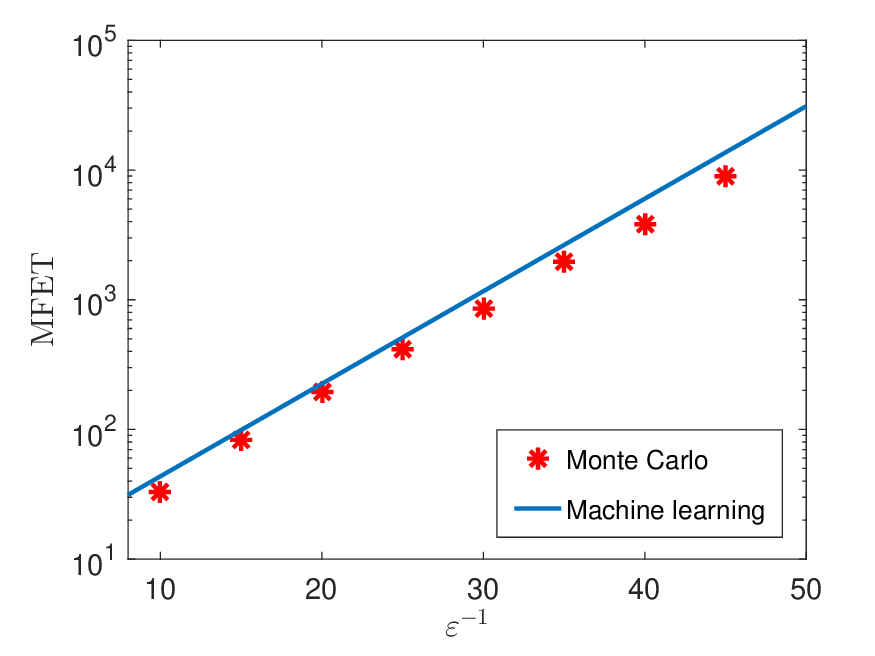}
	\caption{In the characteristic boundary case, the machine learning result, denoted by the blue curve, aligns well with the Monte Carlo simulation  marked by red stars. This concordance reinforces the consistency of both methods and their reliability in predicting the correlation between mean first exit time and noise intensity.}
	\label{fig9}
\end{figure}

Denote $H^{\ast}=\nabla^{2}V_\theta(\text{US})$. We solve the equation
\begin{equation*}
(H^{\ast})^{-1}(Q^{\ast})^{\top}+Q^{\ast}(H^{\ast})^{-1}=a^{\ast},
\end{equation*}
where
\begin{equation*}
Q^{\ast}=-\nabla b(\text{US}),\quad a^{*}=a(\text{US}).
\end{equation*}
We can obtain
\begin{equation*}
H^{\ast}=\left(
           \begin{array}{cc}
           -0.0446   & 0.4228 \\
            0.4228  & 1.3305 \\
           \end{array}
         \right).
\end{equation*}
Then $\det(H^{\ast})=-0.238$. Specifically, $\lambda^{\ast}=0.3721$, $V_{\theta}(\text{US})=0.1643$.
Thus the mean first exit time crossing this boundary can be computed via its asymptotic approximation,
as plotted in Fig. \ref{fig9}. The machine learning results are also validated by Monte Carlo method.
It is observed that the mean first exit time in the characteristic boundary case is much greater than
the one in the non-characteristic boundary case, given the same noise intensity. Therefore, these two boundaries
can be used for two-level early warning.

Finally, we want to investigate how different noise combinations affect the escape behavior of the system.
We consider three cases:
\begin{description}
  \item[(i)] $\sigma_1=0.1$, $\sigma_2=1$;
  \item[(ii)] $\sigma_1=0.08$,  $\sigma_2=0.1$;
  \item[(iii)] $\sigma_1=0.1$, $\sigma_2=0.8$.
\end{description}
We take the characteristic boundary as an illustrative example.
Fig. 10 illustrates the most probable exit paths in the above three cases.
Taking case (i) as a benchmark, we find that the path moves downward
as $\sigma_1$ decreases, while the path moves upward when $\sigma_2$ decreases.
This phenomenon is understandable. When $\sigma_1$ decreases, the effect of noise in the $x_1$ direction becomes smaller, so the path tends to be vertical, allowing the noise in the $x_2$ direction to exert a more significant influence and causing the path to bend downward. Conversely, as $\sigma_2$ decreases, the path becomes more horizontal, causing it to move upward.

\begin{figure}
	\centering
	\includegraphics[width=7cm]{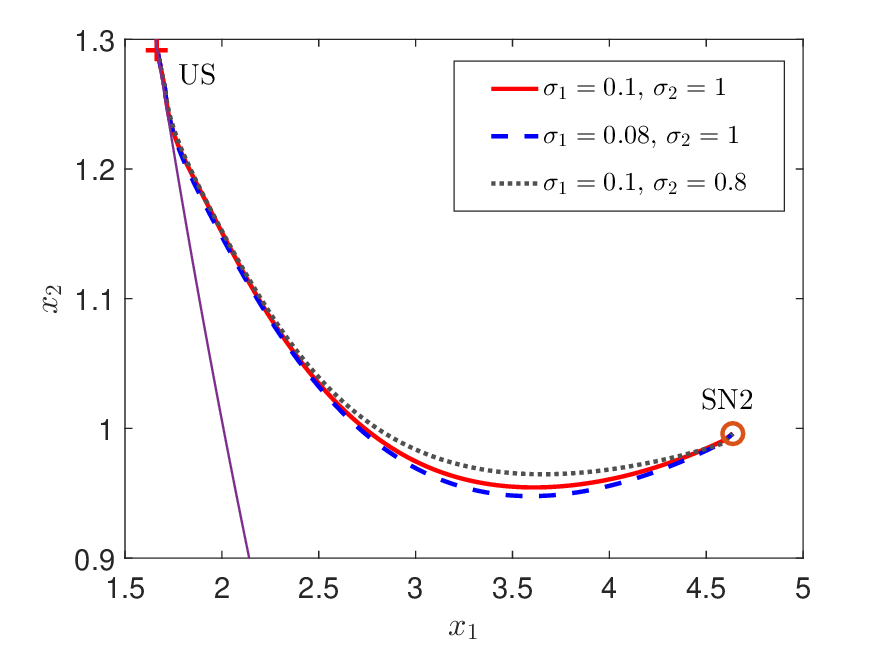}
	\caption{The most probable exit paths from the metastable state SN2 to the saddle point US for the three cases (i)-(iii).}
	\label{fig10}
\end{figure}

\begin{figure}
	\centering
	\includegraphics[width=7cm]{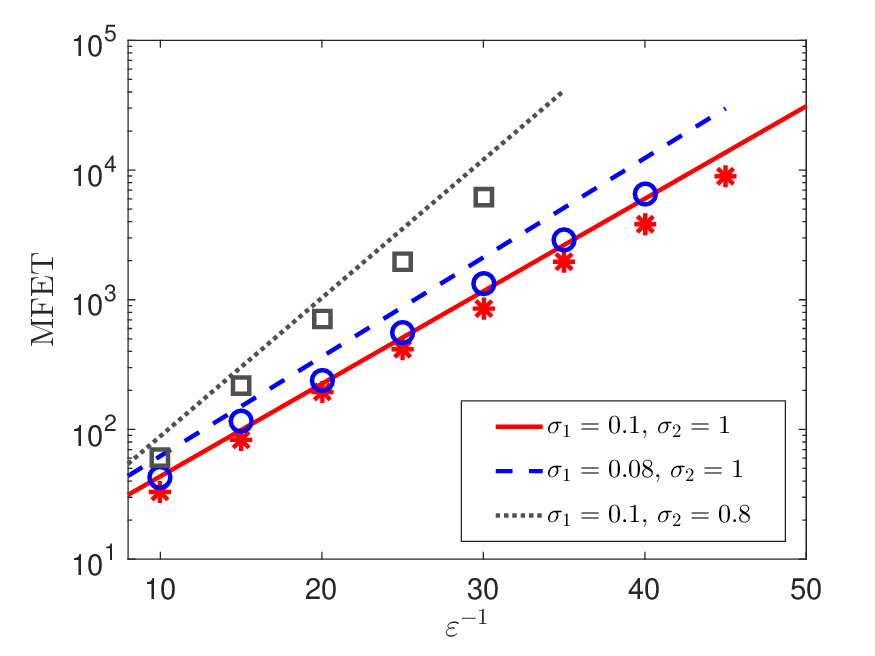}
	\caption{The graph compares the mean first exit time among the three cases (i)-(iii), including a side-by-side  comparison with Monte Carlo simulations.}
	\label{fig11}
\end{figure}

In addition, Fig. \ref{fig11} depicts the mean first exit time of the three cases (i)-(iii), along with a comparison of Monte Carlo simulations.
Obviously, irrespective of the direction in which the noise is reduced, the mean first exit time will increase. Reducing $\sigma_1$ from 0.1 to 0.08 and reducing $\sigma_2$ from 1 to 0.8 result in the same proportional reduction. However, the growth of mean first exit time after reducing $\sigma_2$ is significantly greater than the result of reducing $\sigma_1$. Therefore, the noise in the $x_2$ direction  plays a dominant role in the escape process of the system.

\section{Conclusion and future perspective}\label{C}
In this paper, we proposed a machine learning method to investigate rare events of stochastic dynamical systems with multiplicative Gaussian noise. We computed the most probable paths, quasipotential, and mean first exit time of a stochastic vegetation-water system in both non-characteristic and characteristic boundaries via the machine learning algorithm. We analyzed the dynamics of the system and explored the feasibility of using the exit phenomenon to establish early warnings for vegetation degradation.

We can innovatively apply machine learning methods to more complex and high-dimensional stochastic dynamical systems driven by Gaussian multiplicative noise. Furthermore, rare events, most probable paths, quasipotential, transition rates, and mean first exit time of the dynamical system with non-Gaussian L\'evy noise are also worthy of further study.

\section*{Acknowledgement}
The authors acknowledge support from the National Natural Science Foundation of China (Grant Nos. 12302035, 62073166, 62221004), the Natural Science Foundation of Jiangsu Province (Grant No. BK20220917), the Key Laboratory of Jiangsu Province,  the Shandong Provincial Natural Science Foundation under Grant ZR2021ZD13, and the Project on the Technological Leading Talent Teams Led by Frontiers Science Center for Complex Equipment System Dynamics (FSCCESD220401).

\section*{Data availability}
Numerical algorithms source code
associated with this article can be found, in the online version, at
https://github.com/liyangnuaa/rare-events-in-stochastic-vegetation-system.

\end{document}